\documentclass[twocolumn]{autart}    

\usepackage{graphicx}          
 \usepackage{multirow}                              
\usepackage{amsmath}   		
\usepackage{cases}
\usepackage{graphics} 		
\usepackage{setspace}		
\usepackage{longtable}          
\usepackage{verbatim}
\usepackage{caption}
\usepackage{fancyhdr}      
\usepackage{bm}
\usepackage[noadjust]{cite}
\usepackage{floatrow}
\usepackage{color}
\usepackage{tabularx}

\usepackage{pdfpages}

\usepackage{epsfig} 
\usepackage{mathptmx} 
\usepackage{times} 
\usepackage{amsmath} 
\usepackage{amssymb}  
\usepackage{amsfonts}
\usepackage{tikz,tabstackengine,amsmath}
\usetikzlibrary{shapes,shadows,arrows}
\usetikzlibrary{shapes,shapes.geometric,arrows.meta,fit,calc,positioning,automata}
\usepackage{pgfplots}
\usepgfplotslibrary{groupplots}
\usepackage{subfig}
\usepackage{mathtools}
\usepackage{adjustbox}
\usepackage{xcolor}
\usepackage{color}
\usepackage [english]{babel}
\usepackage [autostyle, english = american]{csquotes}
\usepackage{todonotes}
\usepackage{dsfont}
\newcolumntype{C}{>{\centering\arraybackslash}X}

\usepackage[utf8]{inputenc}

\usepackage{enumitem}

\newcommand{\mb}{\mathbb}

\newcommand{\mc}{\mathcal}

\newcommand{\mbR}{{\mathds{R}}}
\newcommand{\Id}{{\mathbb{I}}}

\newcommand{\hs}[1]{{\hspace{#1cm}}}

\newcommand{\mfN}{{\mathfrak{N}}}

\newcommand{\mfT}{{\mathfrak{T} }}
\newcommand{\mfI}{{\mathfrak{I} }}

\newcommand{\hH}{{\widehat{H}}}

\newcommand{\nn}{{\nonumber}}

\usepackage{stackengine}

\definecolor{myGrn}{rgb}{0,0.5,0}

\usepackage[normalem]{ulem}

\begin{document}

\begin{frontmatter}

\title{LQG Mean Field Games with a Major Agent: Nash Certainty Equivalence versus Probabilistic Approach} 
\thanks[footnoteinfo]{Corresponding author D.~Firoozi.}

\author[Paestum]{Dena Firoozi}\ead{dena.firoozi@hec.ca}    

\address[Paestum]{Department
of Decision Sciences, HEC Montr\'eal, Montreal,
QC, Canada}  

\thanks{The author would like to acknowledge helpful discussions with M. Huang, R. Malham\'{e}, P. E. Caines, and M. Pazoki. }
\begin{abstract} 
Mean field game (MFG) systems consisting of a major agent and a large number of minor agents were introduced in (Huang, 2010) in an LQG setup. The Nash certainty equivalence was used to obtain a Markovian closed-loop Nash equilibrium for the limiting system when the number of minor agents tends to infinity. In the past years several approaches to major--minor mean field game problems have been developed, principally (i) the Nash certainty equivalence and analytic approach, (ii) master equations, (iii) asymptotic solvability, and (iv) the probabilistic approach. For the LQG case, the recent work (Huang, 2021) establishes the equivalency of the Markovian closed-loop Nash equilibrium obtained via (i) with those obtained via (ii) and (iii). In this work, we demonstrate that the Markovian closed-loop Nash equilibrium of (i) is equivalent to that of (iv) for the LQG case. These two studies answer the
long-standing questions about the consistency of the solutions to major-minor LQG MFG systems derived using different approaches.

\end{abstract}
\end{frontmatter}

\section{Introduction}
 Mean field game (MFG) systems with major and minor agents were first introduced in \cite{Huang2010} in an LQG setting, where there is a major agent (whose impact does not vanish in the limit of infinite population size) together with a population of minor agents (where each agent has individually asymptotically negligible effect). In the introduced setting the major agent's state appears in both dynamics and the cost functional of each minor agent. Moreover, each agent is interacting with the average state of minor agents through couplings in the dynamics and cost functionals. As a result, the mean field for such systems is a progressively measurable stochastic process with respect to the filtration generated by the major agent's Wiener process.  
 In \cite{Huang2010}, the author uses the Nash certainty equivalence to establish the existence of Markovian closed-loop $\epsilon$-Nash equilibria and derive the individual agents' explicit control laws which together yield an equilibrium. This methodology is extended in \cite{NourianSiam2013} for a general nonlinear case where
the major agent's state appears in nonlinear dynamics and cost functional of individual minor agents, and all agents are coupled with the empirical distribution of minor agents' state. The best-response strategy of an agent in the limiting case is formulated as the solution to a set of coupled
forward backward (FB) PDEs, i.e.~a Hamilton Jacobi Bellman equation and a Fokker--Planck--Kolmogorov equation. Subsequently, it is shown that the set of best-response strategies yields a Markovian closed-loop $\epsilon$-Nash equilibrium for the system. This methodology which mainly uses the
dynamic programming principle is known as the analytic approach in some literature. 

A probabilistic approach to major-minor (MM) MFG systems by using the stochastic maximum principle is developed in \cite{CarmonaZhu2016}, where the authors establish the existence of open-loop $\epsilon$-Nash equilibria for a general case as the solutions to a set of FBSDEs, and provide the explicit solutions for an LQG case. In \cite[Section 6]{CarmonaZhu2016}, it is discussed in detail that the obtained open-loop equilibrium is different from the Markovian closed-loop equilibrium derived in \cite{Huang2010} for the LQG case. Following this work an alternative probabilistic formulation is proposed in  \cite{CarmonaWang2017}, where the stochastic maximum principle is used and the search for Nash equilibria in the infinite-population limit is formulated as the search for fixed points in the space of best response control maps for the major and minor agents. Using this method, the authors retrieve the same set of FBSDEs as in \cite{CarmonaZhu2016} characterizing the open-loop equilibrium without explicitly solving it. This is while the paper does not present any comparison between the obtained Markovian closed-loop Nash equilibrium and that of the existing work \cite{Huang2010}. Therefore, the paper is inconclusive about this important aspect. 

The MFG master equation methodology, which encapsulates the MFG system in a backward nonlinear PDE in infinite dimension, is used in \cite{LasryLions2018, CardaliaguetCirant2018} to characterize Nash equilibria for a general MM MFG system. 
Moreover, \cite{CardaliaguetCirant2018} shows that the solution of the finite-population MM MFG Nash system converges to the solution of the system of master equations as the number of minor agents tends to infinity. 

The solutions to MM LQG MFG systems obtained using the above discussed methods are seemingly different . Consequently, in order to address the questions in the MFG community related to the consistency of the Nash certainty equivalence solutions \cite{Huang2010}  with the ones  obtained via other methods in the literature, the following works have emerged. \cite{FCJ-Convex2018} uses a variational analysis to retrieve the Markovian solutions of \cite{Huang2010}, where no assumption is imposed on the mean field evolution a priori.  Moreover, \cite{HuangCIS2021} establishes that the Nash certainty equivalence solutions \cite{Huang2010} are equivalent to the Markovian closed-loop solutions obtained via the master equations \cite{LasryLions2018} and asymptotic solvability. The current work serves as the last piece of the long-standing puzzle about MM LQG MFG systems. We demonstrate that the Markovian closed-loop Nash equilibria obtained through the Nash certainty equivalence \cite{Huang2010} and the probabilistic approach \cite{CarmonaWang2017} for the limiting MM LQG MFG systems are equivalent. (For the detailed analysis related to the derivation of consistency equations, best response strategies and $\epsilon$-Nash property we refer the reader to \cite{Huang2010}, \cite{FCJ-Convex2018} and \cite{CarmonaWang2017}.)

In this paper we first introduce finite-population MM LQG MFGs in Section \ref{sec:MMLQGMFG}. Next, we present the Nash certainty equivalence solutions (\cite{FCJ-Convex2018}, \cite{Huang2010}) in Section \ref{sec:NCE}. We then present the Markovian closed-loop solutions obtained via the probabilistic approach (\cite{CarmonaWang2017}) in Section \ref{sec:PA}. Finally, we show that the two solutions are equivalent in Section \ref{sec:comparison}. 

\section{Finite-Population MM LQG MFG Systems}\label{sec:MMLQGMFG}

We consider a large population $N$ of minor agents, each denoted by $\mc{A}_i, i \in \mfN:=\{1,\dots,N\},\,N<\infty$, and a major agent denoted by $\mc{A}_0$. To capture the essence of the two approaches we consider a simple LQG case (for a general case with heterogeneous minor agents see \cite{FCJ-Convex2018, HuangCIS2021,Huang2010}), where the major and minor agents' states, respectively, satisfy
\begin{gather}
dx_t^0 = [A_0\,x_t^0 +F_0\, x^{(N)}_t+ B_0\, u_t^0]\,dt + \sigma_0 \,dw^0_t, \label{MajorDyn}\\
dx_t^i = [A\, x_t^i +F\, x^{(N)}_t + G x^0_t + B\, u_t^i]\,dt + \sigma \,dw^i_t, \label{MinorDyn}
\end{gather}
for $t \in\mfT = [0,T]$, $i \in \mfN$. Here $x^i_t \in \mathbb{R}^n,~ i \in \mfN_0:=\{0,\dots,N\}$, are the states, $(u^i_t)_{t\in\mfT} \in \mathbb{R}^m,~ i \in \mfN_0$, are the control inputs, $w =\lbrace (w^i_t)_{t\in\mfT}, w^i_t \in \mathbb{R}^r, i \in \mfN_0 \rbrace$ denotes $(N+1)$ independent standard Wiener processes. 
 Moreover, $x^{(N)}_t := \tfrac{1}{N} \sum_{i\in \mfN} x^i_t$ denotes the average state of the minor agents. All matrices in \eqref{MajorDyn} and \eqref{MinorDyn} are constant and of appropriate dimension.
\begin{assum} \label{IntialStateAss}
The initial states $\lbrace x^i_0,~ i \in \mfN_0 \rbrace$ defined on $(\Omega, \mathcal{F}, P)$ are identically distributed, mutually independent and also independent of $w$, with $\mathbb{E}x^i_0=\xi,\, i\in \mfN$, and $\sup_{i} \mathbb{E}\Vert x^i_0\Vert^2 \leq c < \infty $, $ i \in \mfN_0$, where $c$ is independent of $N$. Moreover each agent $\mc{A}_i, i \in \mfN_0$, observes $\xi,\, x^0_0$.
\end{assum} 

We denote $\Vert a \Vert_B^2\coloneqq a^\intercal B a$, where $a$ and $B$ are matrices of appropriate dimension. We also denote $u^{-0} \coloneqq  (u^1,\dots,u^N)$ and $u^{-i} \coloneqq  (u^0,\dots,u^{i-1}, u^{i+1},\dots, u^N)$. Then the cost functionals for a major agent and  a minor agent $\mc{A}_i$, $i\in\mfN$, are given by
\begin{gather}
J_{0}^{N}(u^0, u^{-0}) = \tfrac{1}{2} \mathbb{E} \int_{0}^{T}  \Big \lbrace \Vert x^0_t - \Phi^{(N)}_t\Vert ^2_{Q_0} + \Vert u^0_t \Vert_{R_0}^2 \Big \rbrace dt,\label{MajorCostLrgPop}\allowdisplaybreaks\\
J_{i}^{N}(u^i, u^{-i}) = \tfrac{1}{2} \mathbb{E} \int_{0}^{T} \Big \lbrace \Vert x^i_t - \Psi^{(N)}_t \Vert_{Q}^{2} + \Vert u^i_t \Vert_{R}^2 \Big\rbrace \,dt, \label{MinorCostFinitePop}\allowdisplaybreaks\\
\Phi^{(N)}_t := H_0\,x^{(N)}_t + \eta_0, \quad \Psi^{(N)}_t ~:=~ H\, x^0_t + \hH  \, x^{(N)}_t + \eta.
\end{gather}
\begin{assum}(Convexity)\label{ConvexityCondMajorLQGMFG}
$R_0 >0$, $Q_0\geq 0$, $R >0$, $Q \geq 0$.
\end{assum}

\section{Infinite Population MM LQG MFG Systems} The dynamics and cost functional of the major agent $\mc{A}_0$ and a generic minor agent $\mc{A}_i$ in the infinite-population case are given by
\begin{gather}
dx_t^0 = [A_0\,x_t^0 +F_0\, \bar{x}_t+ B_0\, u_t^0]\,dt + \sigma_0 \,dw^0_t, \label{MajorDyn-inf}\allowdisplaybreaks\\
dx_t^i = [A\, x_t^i +F\, \bar{x}_t + G x^0_t + B\, u_t^i]\,dt + \sigma \,dw^i_t, \label{MinorDyn-inf}\allowdisplaybreaks\\
J_{0}^{N}(u^0) = \tfrac{1}{2} \mathbb{E} \int_{0}^{T}  \Big \lbrace \Vert x^0_t - H_0\bar{x}_t + \eta_0\Vert ^2_{Q_0} + \Vert u^0_t \Vert_{R_0}^2 \Big \rbrace dt,\label{MajorCostLrgPop}\allowdisplaybreaks\\
J_{i}^{N}(u^i,u^0) = \tfrac{1}{2} \mathbb{E} \int_{0}^{T} \Big \lbrace \Vert x^i_t -H x^0_t - \hH  \bar{x}_t - \eta \Vert_{Q}^{2} + \Vert u^i_t \Vert_{R}^2 \Big\rbrace \,dt, \label{MinorCostInfPop}
\end{gather}
 where the mean field is defined as $\bar{x}_t := \lim_{N \rightarrow \infty}\frac{1}{N} \sum_{j\in \mfN} x^{j}_t$, if the limit exists. It is equivalently defined as the expected value  $\bar{x}_t=\mathbb{E}[x^i_t|\mathcal{F}^0_t]$ of the state of a generic minor agent $\mc{A}_i$ given the information set $\mathcal{F}^0_t$ defined below.

\text{\textbf{Information Sets.}} We define (i) the major agent's information set $\mc{F}^0:=(\mc{F}^0_t)_{t\in\mfT}$ as the filtration generated by $(w^0_{t})_{t\in\mfT}$, and (ii) a generic minor agent $\mc{A}_i$'s information set $\mc{F}^i\coloneqq (\mc{F}_{t}^i)_{t\in\mfT}$ as the filtration generated by $(w^i_{t}, w^0_t)_{t\in\mfT}$.

\begin{assum}
(Admissible Controls)\label{ass:MajorControl}
 (i) For the major agent $\mc{A}_0$, the set of admissible control inputs  $\mc{U}^{0}$ is defined to be the collection of Markovian linear closed-loop control laws $u^0 \coloneqq (u^0_t)_{t\in\mfT}$ such that $\mb{E}[\int_0^T u_t^{0\intercal}u_t^0\, dt] < \infty$. More specifically, $u^0_t = \ell^0_0(t)+\ell^1_0(t)x^0_t+\ell^2_0(t)\bar{x}_t$ for some deterministic functions $\ell^0_0(t), \ell^1_0(t),$ and $\ell^2_0(t)$. (ii) For each minor agent $\mathcal{A}_i,\, i \in \mfN$, the set of admissible control inputs $\mc{U}^{i}$ is defined to be the collection of Markovian linear closed-loop control laws $u^i\coloneqq (u^i_t)_{t\in\mfT}$ such that  $\mb{E}[\int_0^T u_t^{i\intercal}  u_t^i\, dt] < \infty$. More specifically, $u^i_t = \ell^0(t)+\ell^1(t)x^i_t+\ell^2(t)x^0_t+\ell^3(t)\bar{x}_t$ for some deterministic functions $\ell^0(t),\ell^1(t), \ell^2(t)$ and $\ell^3(t)$.
\end{assum}
\section{Nash Certainty Equivalence Approach}\label{sec:NCE}

In the Nash certainty equivalence approach, first an a priori dynamics for the mean field is derived. Then the idea is to Markovianize (i) the major agent's limiting system by extending its state with the mean field, and (ii) a generic minor agent's limiting system by extending its state with the major agent's state and the mean field. This state extension leads to a set of decoupled classical optimal control problems for individual agents, which are linked with each other through the major agent's state and the mean field. Given the individual information sets, each agent can solve its own stochastic optimal control problem to obtain a best-response strategy. Subsequently, a Nash equilibrium is defined as the set of the best-response Markovian closed-loop strategies of individual agents such that they collectively generate the same mean field that was used in the first step to obtain the best response strategies. This yields a set of consistency equations, the fixed-point solution of which characterizes the Nash equilibrium. (\cite{FCJ-Convex2018} uses a variational analysis and obtains the same Nash equilibrium without assuming an a priori mean field evolution.)

\text{\textbf{Mean Field Evolution.}} According to \cite{Huang2010, HuangCIS2021}, if a generic minor agent adopts a Markovian linear closed-loop strategy,  $\bar{x}_t$ satisfies
\begin{equation}\label{MeanFieldEq}
d\bar{x}_t = \left(\bar{A}\,\bar{x}_t + \bar{G}\, x^{0}_t + \bar{m}(t)\right)dt,
\end{equation}

where $\bar{A},\, \bar{G} \in \mb{R}^{n \times n}$ and $\bar{m}\in\mb{R}^{n}$ are functions of the fixed-point solutions to the consistency equations \eqref{ConsistencyEq1}-\eqref{ConsistencyEq2-4}. Now we present the agents' Markovianized systems.

\textbf{Major Agent.}
From \eqref{MajorDyn-inf} and \eqref{MeanFieldEq} the major agent's extended state $X^0_t = [x^{0\intercal}_t\,\, \bar{x}_t^\intercal ]^\intercal$ satisfies 
\begin{subequations}\label{majorExtDynUnPert}
 \begin{gather}
dX^{0}_t = \left(\mb{A}_0 X^{0}_t + \mb{B}_0 u^0_t  + \mb{M}_0 \right) dt + \Sigma_0 dw^0_t,\\
\mb{A}_0 = \begin{bmatrix}
A_0 & F_0\\
\bar{G} & \bar{A}
\end{bmatrix}, 
\mb{M}_0 = \begin{bmatrix} 0\\ \bar{m} \end{bmatrix},  \mb{B}_0=\begin{bmatrix} B_0 \\ 0  \end{bmatrix}, 
\Sigma_0 = \begin{bmatrix}
\sigma_0 \\
0 
\end{bmatrix}.
\end{gather}
\end{subequations}
The major agent's cost functional in terms of $X^0_t$  is given by
\begin{subequations}
\begin{gather}\label{majorExtCost}
J_0^{\infty}(u^0) = \tfrac{1}{2}\mb{E} \int_0^T \Big \{\Vert X_s^{0}\Vert^2_{\mb{Q}_0}  + \Vert u^0_s \Vert_{R_0}^2 -2(X_s^{0})^\intercal \bar{\eta}_0 \Big \} ds,\allowdisplaybreaks\\
\mathbb{Q}_0 = \Vert\left[\Id_{n}, -H_0\right]\Vert^2_{Q_0}, \quad
\bar{\eta}_0 = \left[\Id_{n}, -H_0\right ]^\intercal  Q_0 \eta_0. 
\end{gather}
\end{subequations}

According to \cite[Thm 5]{FCJ-Convex2018} (the finite-horizon version of \cite[Thm 10]{Huang2010} for general MM LQG MFG systems), the best response strategy for the major agent is given by
\begin{subequations}\label{SextGMFGLatentSLQconvControl}
\begin{gather}
 u^{0,*}_t = - R_0^{-1} \mathbb{B}_0^\intercal \big(\Pi_0(t) X_t^{0,*} + s_0(t) \big),\label{majorCntrl_NCE}\allowdisplaybreaks\\
-\dot{\Pi}_0
= \Pi_0 \mb{A}_0 + \mb{A}_0^\intercal \Pi_0 - \Pi_0 \mb{B}_0R_0^{-1}\mb{B}_0^\intercal \Pi_0+ \mb{Q}_0,\label{SextGMFGlatentLQGriccati}\allowdisplaybreaks\\
-\dot{s}_0 =  \left[\mb{A}_0^\intercal- \Pi_0\, \mb{B}_0\, R_0^{-1} \,\mb{B}_0^\intercal \right] s_0 + \Pi_0 \mb{M}_0 -\bar{\eta}_0,\\
\Pi_0(T) = 0, ~ s_0(T) = 0.
\end{gather}
\end{subequations}

 \textbf{Minor Agent.}
From \eqref{MajorDyn-inf}-\eqref{MinorDyn-inf} and \eqref{MeanFieldEq}, for a generic minor agent $\mc{A}_i$, the extended state $X^{i}_t = [x^{i \intercal}_t, x^{0\intercal}_t, \bar{x}_t^\intercal]^\intercal $ is governed by the dynamcis
\begin{subequations}\label{minorExtDynUnpert}
\begin{gather}
 dX^{i}_t = \left( \mb{A} X^{i}_t  + \mb{B} \,u^i_t  + \mb{M}(t) \right) dt + \Sigma \,dW^i_t,\allowdisplaybreaks\\
 \mb{A} = \begin{bmatrix} A & [G \, \, \, F]\\ 0 &\mb{A}_0-\mb{B}_0 R_0^{-1} \mb{B}_0^\intercal \Pi_0 \end{bmatrix}, ~ \mb{B} = \left[ \begin{array}{c} B \\ 0 \end{array}\right],~ W^i_t = \begin{bmatrix} w^i_t \\ w^0_t \end{bmatrix}, \nonumber\allowdisplaybreaks\\ \mb{M} = \begin{bmatrix} 0 \\ \mb{M}_0-\mb{B}_0 R_0^{-1} \mb{B}_0^\intercal s_0 \end{bmatrix},\,\,
\Sigma = \begin{bmatrix} \sigma & 0 \\ 0 & \Sigma_0 \end{bmatrix}.\label{sysMatMinor}
\end{gather}
\end{subequations}
The cost functional for $\mc{A}_i,\, i\in\mfN$, in terms of $X^i_t$ is given by
\begin{subequations}\label{minorExtCost}
\begin{gather}
J_i^{\infty}(u^i) = \tfrac{1}{2}\mb{E} \int_0^T \Big \{\Vert X_s^{i}\Vert_{\mb{Q}}^2 + \Vert u_s^i \Vert^2_R -2(X_s^{i})^\intercal \,\bar{\eta}
\Big \} ds,\allowdisplaybreaks\\
\mathbb{Q} =  \Vert[\Id_{n},  -H, -\hH]\Vert^2_Q,\quad
\bar{\eta} = [\Id_{n}, -H, -\hH]^\intercal  Q \eta.\label{weightMatMinor}
\end{gather}
\end{subequations}

According to \cite[Thm 5]{FCJ-Convex2018} and  \cite[Thm 10]{Huang2010}, the best response strategy for a generic minor agent $\mc{A}_i$ is given by
\begin{subequations}\label{SextGMFGLatentSLQconvControl-minor}
\begin{gather}
 u^{i,*}_t = - R^{-1} \mathbb{B}^\intercal \big(\Pi(t) \,X_t^{i,*} + s(t) \big),\allowdisplaybreaks\\
-\dot{\Pi} = \Pi \mb{A} + \mb{A}^\intercal \Pi - \Pi \mb{B}R^{-1}\mb{B}^\intercal \Pi + \mb{Q},\quad \Pi(T) = 0,\allowdisplaybreaks\\
- \dot{s}  =  [\mb{A}^\intercal- \Pi \mb{B} R^{-1} \mb{B}^\intercal ] s + \Pi \mathbb{M}- \bar{\eta}, \quad s(T) =  0.
\label{SextGMFGlatentLQGriccati-minor}
\end{gather}
\end{subequations}

\textbf{Mean Field Consistency Equations.}\label{sec:CE} We first define
\begin{equation}
\Pi_k = 
\begin{bmatrix}
\Pi_{11} & \bar{{\Pi}}_{12}\\
\bar{\Pi}_{21} & \bar{\Pi}_{22} 
\end{bmatrix}
, \quad \bar{{\Pi}}_{12}=\begin{bmatrix}{\Pi}_{12} &  {\Pi}_{13}\end{bmatrix}, \quad 
s = 
\begin{bmatrix}
s_{1}\\
\bar{s}_{2}
\end{bmatrix}
,
\end{equation}
where $\Pi_{11}, \Pi_{12}, \Pi_{13} \in \mbR^{n \times n},\, \bar{\Pi}_{22} \in \mbR^{2n \times 2n},\, \bar{\Pi}_{12} \in \mbR^{n\times 2n}$, $\bar{\Pi}_{21} \in \mbR^{2n\times n},\,s_{1}\in \mbR^{n},\,\bar{s}_{2} \in \mbR^{2n}$. Then according to the Nash certainty equivalence approach \cite{Huang2010}, the consistency equations are obtained by effectively equating \eqref{MeanFieldEq} with the mean field equation resulting from the collective action of the mass of minor agents. Subsequently the consistency equations determining $\bar{A}$, $\bar{G}$, $\bar{m}$ are given by 

\begin{numcases}{\hs{-0.9}}\label{ConsistencyEq1}
 -\dot{\Pi}_0  =  \Pi_0 \mb{A}_0 + \mb{A}_0^\intercal \Pi_0 - \Pi_0 \mb{B}_0R_0^{-1}\mb{B}_0^\intercal \Pi_0 + \mb{Q}_0, \label{ConsistencyEq1-1}\\
-\dot{\Pi}  =   \Pi \mb{A} + \mb{A}^\intercal \Pi - \Pi \mb{B}R^{-1}\mb{B}^\intercal \Pi + \mb{Q}, \label{ConsistencyEq1-2}\\
\bar{A}  =  A - B R^{-1} B^\intercal \Pi_{11}  + F - B R^{-1} B^\intercal \Pi_{13}, \\
\bar{G}  = G -B R^{-1}B^\intercal \Pi_{12}, \\
\Pi_0(T) = 0,\quad \Pi(T) = 0,
\end{numcases}

\begin{numcases}{\hs{-0.9}}\label{ConsistencyEq2}
-\dot{s}_0  =   [\mb{A}_0^\intercal- \Pi_0 \mb{B}_0 R_0^{-1} \mb{B}_0^\intercal ] s_0 + \Pi_0\mb{M}_0-\bar{\eta}_0,\label{ConsistencyEq2-1}
\\
-\dot{s}  =  [\mb{A}^\intercal- \Pi \mb{B} R^{-1} \mb{B}^\intercal ] s + \Pi \mathbb{M}- \bar{\eta}, \label{ConsistencyEq2-2}
\\
\bar{m}  =  -B R^{-1} \mathbb{B}^\intercal s,\\
s_0(T) = 0, \quad s(T) = 0.\label{ConsistencyEq2-4}
\end{numcases}

\begin{thm} The consistency equations \eqref{ConsistencyEq1}-\eqref{ConsistencyEq2-4} reduce to 
\begin{numcases}{\hs{-0.9}} 
-\dot{\Pi}_0  = \Pi_0 \mb{A}_0 + \mb{A}_0^\intercal \Pi_0 - \Pi_0 \mb{B}_0R_0^{-1}\mb{B}_0^\intercal \Pi_0 + \mb{Q}_0,\label{CE1maj_red}\allowdisplaybreaks\\
-\dot{\bar{\Pi}}_{12} = \Pi_{11} \left[G\,\, F\right] + A^\intercal \bar{\Pi}_{12} +\bar{\Pi}_{12}  (\mb{A}_0-\mb{B}_0 R_0^{-1} \mb{B}_0^\intercal \Pi_0)\nn\\\hs{1.8} -  \Pi_{11} BR^{-1}B^\intercal \bar{\Pi}_{12}  -Q [H \,\,\, \hH], \label{CE1min_red}\\
\bar{A}  =  A - B R^{-1} B^\intercal \Pi_{11} + F - BR^{-1} B^\intercal \Pi_{13}, \\
\bar{G}  = G -B R^{-1}B^\intercal \Pi_{12},\\
\Pi_0(T) = 0, \,\, \bar{\Pi}_{12}(T) = 0,
\end{numcases}
\begin{numcases}{\hs{-0.9}}
-\dot{s}_0  =   [\mb{A}_0^\intercal- \Pi_0 \mb{B}_0 R_0^{-1} \mb{B}_0^\intercal ] s_0 + \Pi_0\mb{M}_0-\bar{\eta}_0,\label{CE2maj_red}\\
-\dot{s}_{1}  =  [A^\intercal- \Pi_{11} BR^{-1}B^\intercal ] s_{1} - Q\eta \nn\\\hs{2}+ \bar{\Pi}_{12}(\mb{M}_0-\mb{B}_0 R_0^{-1} \mb{B}_0^\intercal s_0),\label{CE2min_red}\\
\bar{m}  =  -B R^{-1} B^\intercal s_{1},\\
s_0(T) = 0, \,\,\, s_{1}(T) = 0,
\end{numcases}
where $\Pi_{11}(T) = 0$ and 
\begin{gather}
-\dot{\Pi}_{11} = \Pi_{11} A + A^\intercal \Pi_{11} -  \Pi_{11} B R^{-1} B^\intercal  \Pi_{11} + Q,\label{CE3_red}\allowdisplaybreaks\\
\bar{\Pi}_{12}=[\Pi_{12},\Pi_{13}], ~~\mb{A}_0 = \begin{bmatrix}
A_0 & F_0\\
\bar{G} & \bar{A}
\end{bmatrix}, 
~~ \mb{B}_0=\begin{bmatrix} B_0 \\ 0  \end{bmatrix},
\nn \allowdisplaybreaks\\
 \mathbb{Q}_0 = \begin{bmatrix} Q_0 & -Q_0H_0\\
-H_0^{\intercal} Q_0 & H_0^{\intercal} Q_0 H_0
 \end{bmatrix}, ~
 \bar{\eta}_0 = \begin{bmatrix}
 Q_0 \eta_0\\
-H_0^{\intercal} Q_0 \eta_0\end{bmatrix},~~ \mb{M}_0 = \begin{bmatrix} 0\\ \bar{m} \end{bmatrix}. \nn
\end{gather}
 \end{thm}
\textit{Proof.}
Given that $\mb{B}^\intercal \Pi = B^\intercal \big[ \Pi_{11}\,\, \bar{\Pi}_{12} \big]$, the optimal control \eqref{SextGMFGLatentSLQconvControl-minor} of the minor agent $\mc{A}_i$ is given by
 \begin{equation}
 {u}^{i,\ast}_t  
= - R^{-1} B^\intercal \Big(
  \Pi_{11} x^{i}_t + \bar{\Pi}_{12} {\big[{x^0_t}^\intercal ~~{\bar{x}_t}^\intercal\big]}^\intercal+ s_{1} \Big). \label{minorCntrlMF_NCE_simp}
 \end{equation}

Hence only the first block row of $\Pi$ and $s$ appear in a generic minor agent's optimal control and the other blocks are irrelevant. Therefore we use \eqref{ConsistencyEq1-2} and \eqref{ConsistencyEq2-2} to derive the equations that $\Pi_{11}$, $\bar{\Pi}_{12}$ and $s_{1}$ satisfy. To this end, we first treat the terms in \eqref{ConsistencyEq1-2} one by one. Block multiplications for the first and the second terms on the right hand side of \eqref{ConsistencyEq1-2} yield 
\begin{gather}
\Pi \mb{A} = 
\begin{bmatrix}
\Pi_{11} A &\,\,\, E_1^\intercal\\
\bar{\Pi}_{21}A &\,\,\, E_2^\intercal
\end{bmatrix}, \vspace{0.2cm}\quad
 \mb{A}^\intercal \Pi = 
\begin{bmatrix}
 A^\intercal \Pi_{11} &\hspace{0.2cm}   A^\intercal \bar{\Pi}_{12}\\
E_1 & E_2
\end{bmatrix},\label{term12Ric}\\
E_1 = [G \, \, F]^\intercal {\Pi}_{11} + (\mb{A}_0-\mb{B}_0 R_0^{-1} \mb{B}_0^\intercal \Pi_0)^\intercal \bar{\Pi}_{21},
\allowdisplaybreaks\\
E_2 = [G \, \, F]^\intercal \bar{\Pi}_{12}+ (\mb{A}_0-\mb{B}_0 R_0^{-1} \mb{B}_0^\intercal \Pi_0)^\intercal \bar{\Pi}_{22}. 
\end{gather} 
For the third and forth terms in \eqref{ConsistencyEq1-2} we have 
\begin{gather}
\Pi \mb{B}R^{-1}\mb{B}^\intercal \Pi =  \begin{bmatrix}
\Pi_{11} BR^{-1}B^\intercal \Pi_{11} & \Pi_{11} BR^{-1}B^\intercal \bar{\Pi}_{12}\\
0 & 0 
\end{bmatrix},\label{term3Ric} \nn\allowdisplaybreaks\\
\mathbb{Q} =\begin{bmatrix}
Q & \hspace{0.2cm} -Q[H\,\, \hH]\\
-[H\,\, \hH]^\intercal Q &\,\,\, [H\,\, \hH]^\intercal Q[H\,\, \hH]
\end{bmatrix}.\label{term4Ric}
\end{gather}
From \eqref{ConsistencyEq1-2} and \eqref{term12Ric}-\eqref{term4Ric}, and through block by block correspondence, we obtain the ODEs that $\Pi_{11}$ and $\bar{\Pi}_{12}$ satisfy as in \eqref{CE3_red} and \eqref{CE1min_red}, respectively.   

Similarly, block multiplications for the terms in \eqref{ConsistencyEq2-2} result in 
\begin{gather}
\Pi \mb{B} R^{-1} \mb{B}^\intercal s = 
\begin{bmatrix}
\Pi_{11} BR^{-1}B^\intercal s_{1}\\
0
\end{bmatrix},\quad \bar{\eta} =
\begin{bmatrix} 
Q\eta\\
-[H \,\,\hH]^\intercal Q\eta
\end{bmatrix}, \nn\allowdisplaybreaks\\
\Pi\mb{M} =
\begin{bmatrix}
\bar{\Pi}_{12}(\mb{M}_0-\mb{B}_0 R_0^{-1} \mb{B}_0^\intercal s_0)\\
\bar{\Pi}_{22} (\mb{M}_0-\mb{B}_0 R_0^{-1}\mb{B}_0^\intercal s_0)
\end{bmatrix}.~
\label{term3Off}
\end{gather}
From \eqref{ConsistencyEq2-2} and \eqref{term3Off}, $s_{1}$ satisfies \eqref{CE2min_red}.
$\hfill \square$

\section{Probabilistic Approach}\label{sec:PA}
In \cite{CarmonaWang2017}, the search for Nash equilibria for MM MFGs is formulated as the search for fixed points in the space of best response control maps for the major and minor agents in the infinite-population limit. In this section, we present the approach of \cite{CarmonaWang2017} for obtaining a Markovian closed-loop equilibrium for MM LQG MFG systems.  

To be self-contained, in Table \ref{table:dyn_procs} we match the notations used in the current work and in \cite{CarmonaWang2017} for presenting the model parameters and the processes. Otherwise, the notations in the two works are the same. 
\begin{table}[ht]\caption{Associated parameters and processes}\vspace{-0.3cm}
\label{table:dyn_procs}
\small
\centering
\begin{tabularx}{1\linewidth}{|p{1.2cm}||p{1.8cm}|p{0.3cm}|p{0.3cm}|p{0.3cm}|p{0.3cm}|p{0.3cm}|p{0.3cm}|}
\cline{1-8}
 Dynamics&Current Work &$A_0$ &$\sigma_0$ & $u^0_t$ & $A$ &$\sigma$ & $u^i_t$ \\[-0.2ex]
  \cline{2-8}
&Reference \cite{CarmonaWang2017}   & $L_0$ & $D_0$  & $\alpha^0_t$ & $L$ & $D$& $\alpha^i_t$\\[-0.2ex]
\cline{1-8}

\end{tabularx}\par\vskip-1.4pt
\begin{tabularx}{\linewidth}{|p{1.2cm}||p{1.8cm}|p{0.45cm}|p{0.45cm}|p{0.45cm}|p{0.45cm}|p{0.45cm}|}
\cline{1-7}
 \multirow{2}{*}{Cost} &Current Work &$Q_0$ &$Q$ & $R_0$ & $R$ &$\hat{H}$\\[-0.2ex]
 \cline{2-7}
 &Reference \cite{CarmonaWang2017}   & $2Q_0$ & $2Q$  & $2R_0$ & $2R$ & $H_1$\\[-0.2ex]
  \cline{1-7}
\end{tabularx}
\end{table}
In \cite{CarmonaWang2017} Markovian linear closed-loop control actions as in Assumption \ref{ass:MajorControl} are considered for the major agent and a representative minor agent, denoted, respectively, by $\alpha^0_t$ and $\alpha^i_t$. The mean field dynamics is then obtained by forming the closed-loop system for the representative agent $\mathcal{A}_i$ using $\alpha^i_t$ and taking the conditional expectation $\mathbb{E}[x^i_t|\mathcal{F}^0_t]$ of its state $x^i_t$ given the information set $\mathcal{F}^0_t$. Subsequently, to obtain the solutions to the major agent's problem, its state is extended with the mean field in the same manner as in \cite{Huang2010}. Then using the stochastic maximum principle for the extended system, the major agent's optimal control action is obtained as 
\begin{equation}\alpha^{0,\ast}_t = - (2R_0)^{-1}\big[0\,\, B_0^\intercal \big] \mb{Y}_t,\label{majorCntrl_PA}\end{equation}
and a set of McKean-Vlasov FBSDEs is derived which solves for the major agent's extended state and the decoupling field (adjoint process) $\mb{Y}_t$. To solve the FBSDEs, an ansatz is adopted for $\mb{Y}_t$ as in  
\begin{equation}\label{majorAdj_PA}
\mb{Y}_t = K_t {\big[{\bar{x}_t}^\intercal ~~{x^0_t}^\intercal\big]}^\intercal+ k_t, 
\end{equation}
where $K_t, k_t,$ are deterministic matrices of appropriate dimension. Then a set of ODEs that $K_t$ and $k_t$ satisfy is derived. Subsequently the notion of a deviating minor agent is introduced as an extra virtual minor agent who deviates from the strategy of its peers and aims to optimize in response to the major agent and the rest of minor agents. However, when dealing with the optimal control problem for the deviating minor agent, the major agent's state and the mean field are considered as exogenous stochastic coefficients, which are determined offline by solving a set of SDEs, i.e. the major agent's closed-loop extended system. Then using the stochastic maximum principle, the deviating minor agent's optimal control is obtained as in
 \begin{equation}\label{minorCntrl_PA}\alpha^{i,\ast}_t = (2R)^{-1} B^\intercal Y^i_t,\end{equation}
 and a set of FBSDEs with random coefficients (not of McKean-Vlasov type), that the minor agent's state $x^i_t$ and decoupling field (adjoint process) $Y^i_t$ satisfy, are derived. To solve the FBSDEs, an ansatz is considered for $Y_t^i$ as in  
\begin{equation}\label{minorAdj_PA}
Y^i_t = \mb{S}_t  {\big[{\bar{x}_t}^\intercal ~~{x^0_t}^\intercal\big]}^\intercal + S_t x^i_t + \bar{s}_t, 
\end{equation}
where $\mb{S}_t$, $S_t$, and $\bar{s}_t$ are matrices of appropriate dimension. The mean field equation resulting from a minor agent using $\alpha^{i,\ast}_t$ must match the one obtained using $\alpha^i_t$ to solve for the exogenous stochastic coefficients $[{\bar{x}_t}^\intercal ~ {x^0_t}^\intercal ]^\intercal$ in the minor agent's optimal control problem. Subsequently the consistency equations whose fixed-point solutions determine $K_t, k_t$ and $\mb{S}_t, S_t, \bar{s}_t$, are given by  (\cite[eq. (31)-(32)]{CarmonaWang2017})
\begin{numcases}{\hs{-0.9}}
-\dot{K}_t = K_t [\mb{L}_t - \bar{\mb{B}}(2R)^{-1}B^\intercal \mb{S}_t] - K_t \bar{\mb{B}}_0(2R_0)^{-1}\bar{\mb{B}}^\intercal K_t  \nn\\ \hs{1}+[\mb{L}_t - \bar{\mb{B}} (2R)^{-1} B^{\intercal} \mb{S}_t]^{\intercal}K_t + 2 \mb{F}_0,\label{CE1maj_CW}\label{CE1maj_CW}\\
-\dot{\mb{S}}_t = \mb{S}_t[\mb{L}_t - \bar{\mb{B}}(2R)^{-1}B^\intercal \mb{S}_t] - \mb{S}_t \bar{\mb{B}}_0 (2R_0)^{-1} \bar{\mb{B}}_0^\intercal K_t\nn\\\hs{1}+ [L^\intercal - S_t B (2R)^{-1}B^\intercal] \mb{S}_t  \nn\\\hs{1}+ [S_t F -2 Q H_1\,\,\, S_t G -2 Q H],\label{CE1min_CW}\\
K_T =0, \,\, \mb{S}_T =0,
\end{numcases}

\begin{numcases}{\hs{-0.9}}
-\dot{k}_t = [\mb{L}_t - \bar{\mb{B}}(2R)^{-1}B^\intercal \mb{S}_t]k_t - K_t \bar{\mb{B}}_0 (2R_0)^{-1} \bar{\mb{B}}_0^\intercal k_t \nn\\\hspace{1cm}- K_t \bar{\mb{B}} (2R)^{-1} B^\intercal \bar{s}_t + 2f_0,\label{CE2maj_CW}\\
-\dot{\bar{s}}_t = [L^\intercal - S_t B (2R)^{-1} B^\intercal]\bar{s}_t- \mb{S}_t \bar{\mb{B}}_0(2R_0)^{-1}\bar{\mb{B}}_0^\intercal k_t \nn\\\hspace{1cm}- \mb{S}_t \bar{\mb{B}}_0(2R)^{-1} B^\intercal \bar{s}_t - 2 Q \eta,\label{CE2min_CW}\\
k_T =0, \,\, \bar{s}_T=0,
\end{numcases}
where 
\begin{gather}
-\dot{S}_t = S_t L + L^\intercal S_t - S_t B (2R)^{-1} B^{\intercal} S_t + 2Q, \quad S_T =0, \label{CE3_CW}\\
\mb{L} = \begin{bmatrix}
L+F-B(2R)^{-1}B^\intercal S_t & G\\
F_0 & L_0 
\end{bmatrix},
 \,\,\, 
\bar{\mb{B}} = \begin{bmatrix}
B\\
0
\end{bmatrix},\,\,\, \bar{\mb{B}}_0 = \begin{bmatrix}
0\\
B_0
\end{bmatrix},\,\, \nn\allowdisplaybreaks\\
\mb{F}_0 = \begin{bmatrix}
H_0^\intercal Q_0 H_0 & -H_0^\intercal Q_0\\
-Q_0H_0 & Q_0
\end{bmatrix}, \,\, ~
f_0 = \begin{bmatrix}
H_0^\intercal Q_0 \eta_0\\
-Q_0 \eta_0
\end{bmatrix}.
\end{gather}
\section{Comparison of the Two Approaches}\label{sec:comparison}
We start with the following theorem.
\begin{thm}\label{thm:equivalency} For the MM LQG MFG system \eqref{MajorDyn-inf}-\eqref{MinorCostInfPop}, the Markovian closed-loop Nash equilibrium obtained via the Nash certainty equivalence is  equivalent to the one obtained via the probabilistic approach. \end{thm}

\textit{Proof.} By inspection, the Markovian linear optimal control laws $\{u^{0,\ast}, u^{i,\ast}, i\in \mfN\}$ obtained through the Nash certainty equivalence (see \eqref{majorCntrl_NCE}, \eqref{minorCntrlMF_NCE_simp}) have the same structure as the ones $\{\alpha^{0,\ast}, \alpha^{i,\ast}, i\in \mfN\}$ obtained through the probabilistic approach (see \eqref{majorCntrl_PA}-\eqref{minorAdj_PA}). It remains to show the equivalency of the sets of consistency equations, the fixed-point solutions of which yield the coefficients in the above control laws. More specifically, we show that the reduced consistency equations \eqref{CE1maj_red}-\eqref{CE3_red} obtained via the Nash certainty equivalence are the same as the consistency equations  \eqref{CE1maj_CW}-\eqref{CE3_CW} obtained via the probabilistic approach. 
To this end, we first define a block elementary operator which operates on a matrix to produce the desired interchanged block rows, as in 
\begin{equation}
\mfI = \begin{bmatrix} 0 & \Id\\
\Id & 0
\end{bmatrix},
\end{equation}
where the identity matrices $\mb{I}$ are of appropriate dimension. Then we correspond the processes in \eqref{CE1maj_red}-\eqref{CE3_red} with those in   \eqref{CE1maj_CW}-\eqref{CE3_CW} as shown in Table \ref{table:CE_procs}.
\begin{table}[ht]
\caption{Corresponding processes in consistency equations}\vspace{-0.3cm}
\begin{tabular}{|p{2cm}||p{1.1cm}|p{0.65cm}|p{.9cm}|p{0.6cm}|p{0.4cm}|}
 \hline
\hspace{-0.15cm}Current Work &$\!\!\!\mfI^\intercal\Pi_0(t)\mfI$ &$ \!\!\!\Pi_{11}(t)$ & $ \!\!\!\bar{\Pi}_{12}(t) \mfI$ &$ \!\!\!s_0(t)\mfI$ & $ \!\!\!s_1(t)$\\[-0.2ex]
 \hline
\hspace{-0.25cm} Reference \cite{CarmonaWang2017}   & $K_t$  & $S_t$ & $\mb{S}_t$& $k_t$ & $\bar{s}_t$ \\[-0.2ex]
 \hline
\end{tabular}
\label{table:CE_procs}
\end{table}
Using Tables \ref{table:dyn_procs}-\ref{table:CE_procs}, we can retrieve  \eqref{CE3_CW} from \eqref{CE3_red} or vice versa. Now we correspond the terms in \eqref{CE1maj_red}-\eqref{CE1min_red} to those in \eqref{CE1maj_CW}-\eqref{CE1min_CW}. First we use Table \ref{table:dyn_procs} to replace the system parameters in \eqref{CE1maj_red} with those considered in \cite{CarmonaWang2017} (i.e. replace $Q_0, R_0,$ respectively, with $2Q_0, 2R_0,$ in \eqref{CE1maj_red}). Then as per Table \ref{table:CE_procs} to obtain $K_t$ we multiply both sides of \eqref{CE1maj_red} from right by $\mfI$ and from left by $\mfI^\intercal$ as in 
\begin{equation*}
- \mfI^\intercal \dot{\Pi}_0 \mfI  =  \mfI^\intercal (\Pi_0 \mb{A}_0 +  \mb{A}_0^\intercal \Pi_0  -  \Pi_0 \mb{B}_0(2R_0)^{-1}\mb{B}_0^\intercal \Pi_0 + 2 \mb{Q}_0) \mfI,
\end{equation*}  
which by inspection is the same equation as \eqref{CE1maj_CW}, particularly 
 \begin{gather}
\mfI^\intercal \Pi_0 \mb{A}_0 \mfI = K_t [\mb{L}_t - \bar{\mb{B}}(2R)^{-1}B^\intercal \mb{S}_t],\quad 2\mfI^\intercal \mb{Q}_0 \mfI = 2\mb{F}_0, 
\nn\allowdisplaybreaks\\
\mfI^\intercal \Pi_0 \mb{B}_0(2R_0)^{-1}\mb{B}_0^\intercal \Pi_0 \mfI = K_t \bar{\mb{B}}_0(2R_0)^{-1}\bar{\mb{B}}^\intercal K_t. \end{gather}
 Next we use Tables \ref{table:dyn_procs}-\ref{table:CE_procs} to match \eqref{CE1min_red} and \eqref{CE1min_CW}. We first replace $R_0, R, Q, \hat{H}$, respectively with $2R_0, 2R, 2Q, H_1$ in \eqref{CE1min_red}, and then right multiply both sides by $\mfI$ which gives us equation \eqref{CE1min_CW}.
Now we show that \eqref{CE2maj_CW} can be retrieved from \eqref{CE2maj_red}. From Tables \ref{table:dyn_procs}-\ref{table:CE_procs}, we first replace $R_0, Q_0$ with $2R_0, 2Q_0$ in \eqref{CE2maj_red}, and then left multiply its both sides by $\mfI$ as in
\begin{equation} 
- \mfI \dot{s}_0 =    \mfI ([\mb{A}_0^\intercal- \Pi_0 \mb{B}_0 (2R_0)^{-1} \mb{B}_0^\intercal ]  s_0 + \Pi_0\mb{M}_0 - 2 \bar{\eta}_0), 
\end{equation}
where $-2 \mfI \bar{\eta}_0 = 2f_0$ and 
\begin{equation*}
\mfI \Pi_0\mb{M}_0 = -K_t \begin{bmatrix}B(2R)^{-1} B^\intercal \bar{s}_t\\0 \end{bmatrix}
= - K_t \bar{\mb{B}}(2R)^{-1} B^\intercal \bar{s}_t. \\ 
\end{equation*}
Finally, we replace $R_0, R, Q,$ respectively, with $2R_0, 2R, 2Q$, in \eqref{CE2min_red}, which yields \eqref{CE2min_CW}. This matches \eqref{CE2min_red} and \eqref{CE2min_CW}. $\hfill \square$

\textbf{Discussions.} To obtain Markovian closed-loop Nash equilibria for MM LQG MFG systems, both Nash certainty equivalence and probabilistic approaches assume that a generic minor agent adopts a Markovian linear closed-loop strategy. With this assumption the mean field equation is derived using an ansatz for the minor agent's control action. Then to solve the major agent's limiting optimal control problem, its state is extended with the mean field in both approaches. In \cite{Huang2010}, the optimal linear state feedback control for the major agent's extended system is obtained using the known results for single-agent LQG systems. This is while in \cite{CarmonaWang2017}, the stochastic maximum principle is used, where a linear ansatz for the adjoint process in terms of the major agent's state and the mean field is considered. The two methods are equivalent and result in the same optimal control for the major agent. For a generic minor agent's optimal control problem, \cite{Huang2010} Markovianizes the minor agent's system by extending its state by the major agent's state and the mean field. Then again using the known results for single-agent LQG systems the minor agent's optimal control is obtained, which is a linear function of its own state, the major agent's state and the mean field. This is while in \cite{CarmonaWang2017}, the major agent's state and the mean field are considered as exogenous stochastic coefficients in the minor agent's system. Then using the stochastic maximum principle, an optimal control is obtained for the minor agent by adopting an ansatz for the adjoint process which is a linear function of its own state, the major agent's state and the mean field. Although, the obtained optimal control actions for a generic minor agent derived in \cite{Huang2010} and \cite{CarmonaWang2017} do not look the same at first glance, \textit{Theorem \ref{thm:equivalency}} establishes that they are indeed equivalent. Hence both approaches yield the same Markovian closed-loop Nash equilibrium for MM LQG MFGs.

The fact that the set of consistency equations of \cite{Huang2010} reduces to that of \cite{CarmonaWang2017} stems from an interaction asymmetry in the minor agent's extended system in the former. In fact, in the minor agent's extended system, the individual minor agent's state and control action do not affect the joint system of the major agent and the mean field (i.e. the major agent's extended system). However, the major agent's state and the mean field affect the dynamics and the cost functional of the individual minor agent. In the core, the minor agent's extended system (modelled in \cite{Huang2010}) is working in the same manner as the individual minor agent system with exogenous stochastic coefficients solving the major agent's extended system (modelled in \cite{CarmonaWang2017}). Such asymmetric interactions do not occur in the major agent's extended system. This is due to the mutual interactions as the major agent's state appears in the mean field dynamics and the mean field appears in both the major agent's dynamics and cost functional.

\bibliographystyle{plain}        
\bibliography{bib_Dena_21Oct20}           

\begin{thebibliography}{1}

\bibitem{CardaliaguetCirant2018}
P.~Cardaliaguet, M.~Cirant, and A.~Porretta.
\newblock Remarks on {Nash} equilibria in mean field game models with a major
  player.
\newblock {\em Proceedings of the American Mathematical Society},
  148(10):4241--4255, 2020.

\bibitem{CarmonaWang2017}
R.~Carmona and P.~Wang.
\newblock An alternative approach to mean field game with major and minor
  players, and applications to herders impacts.
\newblock {\em Applied Mathematics \& Optimization}, 76(1):5--27, 2017.

\bibitem{CarmonaZhu2016}
R.~Carmona and X.~Zhu.
\newblock A probabilistic approach to mean field games with major and minor
  players.
\newblock {\em Annals of Applied Probability}, 26(3):1535--1580, 2016.

\bibitem{FCJ-Convex2018}
D.~Firoozi, S.~Jaimungal, and P.~E. Caines.
\newblock Convex analysis for {LQG} systems with applications to major--minor
  {LQG} mean--field game systems.
\newblock {\em Systems \& Control Letters}, 142:104734, 2020.

\bibitem{Huang2010}
M.~Huang.
\newblock Large-population {LQG} games involving a major player: The {N}ash
  certainty equivalence principle.
\newblock {\em SIAM Journal on Control and Optimization}, 48(5):3318--3353,
  2010.

\bibitem{HuangCIS2021}
M.~Huang.
\newblock Linear-quadratic mean field games with a major player: Nash certainty
  equivalence versus master equations.
\newblock {\em Communications in Information and Systems}, 21(3):441--471,
  2021.

\bibitem{LasryLions2018}
J.~M. Lasry and P.~L. Lions.
\newblock Mean-field games with a major player.
\newblock {\em Comptes Rendus Mathematique}, 356(8):886 -- 890, 2018.

\bibitem{NourianSiam2013}
M.~Nourian and P.~E. Caines.
\newblock $\epsilon$-{N}ash mean field game theory for nonlinear stochastic
  dynamical systems with major and minor agents.
\newblock {\em SIAM Journal on Control and Optimization}, 51(4):3302--3331,
  2013.

\end{thebibliography}

\end{document}